\documentclass[11pt]{article}
\usepackage{amscd,amsmath, amssymb, fancyhdr, color}

\newcommand{\version}{version 2,\ \ October 6, 2014}

\newcommand{\la}{\lambda}

\newcommand{\ka}{K\"ahler}

\newcommand{\te}{\theta}

\newcommand{\mi}{\mathrm{i}}

\newcommand{\norm}[1]{\Vert #1\Vert}

\newcommand{\vv}{\overrightarrow}
\newcommand{\Ker}{\text{Ker}}

\numberwithin{equation}{section}

\def\eqref#1{(\ref{#1})}

\newcommand{\C}{{\mathbb C}}
\newcommand{\R}{{\mathbb R}}

\newcommand{\6}{\partial}
\def\1{\sqrt{-1}\:}

\newcommand{\cntrct}                % contraction with a vector field
{\hspace{2pt}\raisebox{1pt}{\text{$\lrcorner$}}\hspace{2pt}}

\newcommand{\calo}{{\cal O}}

\newcommand{\calh}{{\cal H}}
\newcommand{\cale}{{\cal E}}
\newcommand{\calf}{{\cal F}}

\renewcommand{\bar}{\overline}
\renewcommand{\phi}{\varphi}
\renewcommand{\epsilon}{\varepsilon}
\renewcommand{\geq}{\geqslant}
\renewcommand{\leq}{\leqslant}

\renewcommand{\Re}{\operatorname{Re}}
\renewcommand{\Im}{\operatorname{Im}}

\newcounter{Mycounter}[section]
\newcounter{lemma}[section]
\renewcommand{\thelemma}{{Lemma \thesection.\arabic{lemma}}}
\newcommand{\lemma}{%
     \setcounter{lemma}{\value{Mycounter}}
     \refstepcounter{lemma}
     \stepcounter{Mycounter}
     {\noindent \bf \thelemma:\ }}

\newcounter{claim}[section]
\renewcommand{\theclaim}{{Claim \thesection.\arabic{claim}}}
\newcommand{\claim}{%
     \setcounter{claim}{\value{Mycounter}}
     \refstepcounter{claim}
     \stepcounter{Mycounter}
     {\noindent \bf \theclaim:\ }}

\newcounter{sublemma}[section]
\newcounter{corollary}[section]
\renewcommand{\thecorollary}{{Corollary \thesection.\arabic{corollary}}}
\newcommand{\corollary}{%
     \setcounter{corollary}{\value{Mycounter}}
     \refstepcounter{corollary}
     \stepcounter{Mycounter}
     {\noindent \bf \thecorollary:\ }}

\newcounter{theorem}[section]
\renewcommand{\thetheorem}{{Theorem \thesection.\arabic{theorem}}}
\newcommand{\theorem}{%
     \setcounter{theorem}{\value{Mycounter}}
     \refstepcounter{theorem}
     \stepcounter{Mycounter}
     {\noindent \bf \thetheorem:\ }}

\newcounter{conjecture}[section]
\newcounter{proposition}[section]
\renewcommand{\theproposition}
       {{Proposition \thesection.\arabic{proposition}}}
\newcommand{\proposition}{%
     \setcounter{proposition}{\value{Mycounter}}
     \refstepcounter{proposition}
     \stepcounter{Mycounter}
     {\noindent \bf \theproposition:\ }}

\newcounter{definition}[section]
\renewcommand{\thedefinition}
       {{Definition~\thesection.\arabic{definition}}}
\newcommand{\definition}{%
     \setcounter{definition}{\value{Mycounter}}
     \refstepcounter{definition}
     \stepcounter{Mycounter}
     {\noindent \bf \thedefinition:\ }}

\newcounter{example}[section]
\newcounter{remark}[section]
\renewcommand{\theremark}{{Remark \thesection.\arabic{remark}}}
\newcommand{\remark}{%
     \setcounter{remark}{\value{Mycounter}}
     \refstepcounter{remark}
     \stepcounter{Mycounter}
     {\noindent \bf \theremark:\ }}

\newcounter{problem}[section]
\newcounter{question}[section]
\makeatletter

\@addtoreset{equation}{section}
\@addtoreset{footnote}{section}
\makeatother

\def\blacksquare{\hbox{\vrule width 5pt height 5pt depth 0pt}}
\def\endproof{\blacksquare}

\begin{document}

\begin{center}
{\LARGE\bf
Currents on locally conformally K\"ahler manifolds}\\[3mm]
%%%%%%%%%%%%%%%%%%%%%%%%%%%%%%%%%%%%%%%%%%%%%%%%%%%%%%%%%%%%
{\Large 
Alexandra Otiman\footnote{Partially supported by CNCS grant RU-TE-2011-3-0053.}}
\\[5mm]
\end{center}



\hfill

%%%%%%%%%%%%%%%%%%%%%%%%%%%%%%%%%%%%%%%%%%%%%%%%
{\small
%\hspace{0.15\linewidth}
\noindent\begin{minipage}[t]{0.7\linewidth}
{\bf Abstract} We characterize the existence of a locally conformally \ka\ metric on a compact complex manifold in terms of currents, adapting the celebrated result of Harvey and Lawson for \ka\ metrics.\\ 
\end{minipage}
}
%%%%%%%%%%%%%%%%%%%%%%%%%%%%%%%%%%%%%%%%%%%%%%%%
%\tableofcontents
%%%%%%%%%%%%%%%%%%%%%%%%%%%%%%%%%%%%%%%%%%%%%%%%
%\section{Introduction}
%%%%%%%%%%%%%%%%%%%%%%%%%%%%%%%%%%%%%%%%%%%%%%%%

%%%%%%%%%%%%%%%%%%%%%%%%%%%%%%%%%%%%%%%%%%%%%%%%%%%%%%%%%%%%%%%%%%%%%%%%
%\subsection{Bimeromorphic maps and locally conformally K\"ahler structures}
%%%%%%%%%%%%%%%%%%%%%%%%%%%%%%%%%%%%%%%%%%%%%%%%%%%%%%%%%%%%%%%%%%%%%%%%
\hfill

\hfill
\section{Introduction}

A locally conformally \ka\ manifold (LCK for short) is a Hermitian manifold $(M,J,g)$ for which the fundamental two-form $\omega(X,Y)=g(JX, Y)$ satisfies
\begin{equation}\label{LCK}
d\omega=\theta\wedge\omega,\qquad d\theta=0
\end{equation}
for some one-form $\theta$ called the Lee form. 

There are many examples of compact LCK and non-\ka\ manifolds, among them the Hopf manifolds, see \cite{do}, \cite{ov_survey}.

As $d\theta=0$, the twisted differential  $d_\theta:=d-\theta\wedge$ defines a twisted cohomology which is the Morse-Novikov cohomology of $X$. The LCK condition simply means that the fundamental form of $(X,J,g)$ is $d_\theta$-closed. 

The aim of this note is to obtain an analogue of the intrinsic characterization in \cite{hl} for \ka\ manifolds in the context of LCK geometry. 

\section{LCK condition in terms of currents}

Our main result is the following:

\hfill

\theorem \label{main}
{\em Let $X$ be a compact, complex manifold of complex dimension $n\geq 2$, and let $\theta$ be a closed one-form on $X$. Then $X$ admits a LCK metric with Lee form $\theta$ if and only if there are no non-trivial positive currents which are $(1,1)$ components of $d_\theta$-boundaries.}

\hfill

\remark{Suppose X is a compact complex manifold, admitting a LCK metric, $\omega$, with Lee form $\te$. Then any closed 1-form $\eta \in [\te]_{dR}$ will be a Lee form for a conformal metric 
of $\omega$ and moreover, any conformal change of $\omega$ will be LCK with a Lee form in the same de Rham cohomology class as $\te$.  Therefore, we need not fix $\te$, we can directly use its  
cohomology class, $[\te]_{dR}$. By this observation, the theorem above can be stated as:   

\em Let $X$ be a compact, complex manifold of complex dimension $n\geq 2$, and let $[\theta]_{dR}$ a cohomology class in $H^1_{dR}(X)$. 
Then $X$ admits a LCK metric with Lee form $\theta$ if and only if there are no non-trivial positive currents which are $(1,1)$-components of $d_\eta$-boundaries, 
for any closed one-form $\eta$ belonging to $[\te]_{dR}$. }

\hfill

The rest of Section 2 is devoted to the proof, which follows the lines in \cite{hl}. We use the same results and intermediate steps as \cite{hl}, 
the difficult part being that of finding  some proper analogues in LCK geometry for the \ka\ notions used in the original article. Each following subsection is a step of the proof.

\subsection{Range of $d_\te$ is closed}

Associated with $d_\theta$ are the following operators:
$$\6_\te=\6-\te^{1,0}\wedge,\qquad \bar\6_\te=\bar\6-\te^{0,1}\wedge,\qquad d_\te^c=\mi (\6_\te-\bar\6_\te)$$

\hfill

\definition
A smooth function is called $\theta$-pluriharmonic if it is locally the real part of a smooth $\bar\6_\te$-closed function.

\hfill

We let $\calh_\te$ be the sheaf of germs of $\te$-pluriharmonic functions on $X$.

\hfill

\lemma\label{calh}
{\em $\calh_\te$ is the kernel of the sheaves morphism $\displaystyle \cale_\R\stackrel{d_\te d_{\te}^c}{\xrightarrow{\hspace*{0,5 cm}}}\cale^{1,1}_\R$, where the subscript $\R$ denotes the germs of real valued forms.}

\noindent{\bf Proof:} 
The proof is based on the following easy observation
$$\bar\6_\te f=0 \Leftrightarrow\frac{1}{2\mi}\big(\bar\6_\te f-\6_\te \bar f\big)=0$$
Let now $f=u+\mi v$. One obviously has
\begin{equation}\label{*}
\bar\6_\te f=0 \,\Leftrightarrow\, \frac{\bar\6_\te(u+\mi v)-\6_\te (u-\mi v)}{2\mi}=0 \,\Leftrightarrow\, d_\te v+d_\te ^cu=0
\end{equation}

Let $g$ for which a $f'$ exists such that  $f=g+\mi f'$ is $\bar\6_\te$-closed. It follows from \eqref{*} that $d_\te f'+ d_\te ^cg=0$, which implies $d_\te d_\te ^cg=0$.

Conversely, if $g$ satisfies $d_\te d_\te ^cg=0$, finding a $f'$ such that $f=g+\mi f'$ is $\bar\6_\te$-closed is equivalent to solving the equation $d_\te f'=-d_\te ^c g$.

Since $\te$ is locally exact, let $\te=dh$ on a contractible open set. Then $e^{-h}d_\te ^cg$ is closed and by Poincar\'e lemma there exists a function $h'$ such that $e^{-h}d_\te ^cg =dh'$. Then $f'=-e^hh'$ which completes the proof.

\hfill \endproof

The above result shows that the following is an exact sequence of sheaves:

\begin{equation*}
0\longrightarrow \calh_\te\longrightarrow \cale_\R\stackrel{d_\te d_{\te}^c}{\xrightarrow{\hspace*{1cm}}}\cale^{1,1}_\R \stackrel{d_\te}{\longrightarrow} \big[\cale^{1,2}\oplus\cale^{2,1}\big]_\R\stackrel{d_\te}{\longrightarrow}\cdots
\end{equation*}

Since $\big[\cale^{p,q}\big]_\R$ are acyclic, the above is a resolution which computes the cohomology groups of $\calh_\te$.

\hfill

We now prove that $H^i(X,\calh_\te)$ are finite dimensional for all $i\geq 0$.

Let $\calo_\te$ denote the sheaf of germs of smooth functions satisfying $\bar\6_{\te} f=0$ and let $\calf$ be the kernel of the sheaves morphism $\Re:\calo_\te\rightarrow\calh_\te$:

\begin{equation}\label{**}
0\longrightarrow \calf \longrightarrow \calo_\te\stackrel{\Re}{\longrightarrow}\calh_\te \longrightarrow 0
\end{equation}

\proposition
{\em $\calo_\te$ is locally free of rank $1$ over the sheaf of germs of holomorphic functions, $\mathcal{O}_{X}$, and $\calf$ is locally constant.}

\noindent{\bf Proof:} To prove that $\calf$ is locally constant, we characterize the non-zero $\bar\6_{\theta}$-closed real valued functions. 

Let $h$ be a (unique up to addition with constants) real valued smooth function on a contractible neighbourhood such that $\te^{0,1}=\bar\6 h$. Then
$$\bar\6_{\te}f=0\Leftrightarrow \bar\6 f-f\te^{0,1}=0\Leftrightarrow \bar\6 f=f\bar\6 h$$
Since both $f$ and $h$ are real valued, the above last equality gives, by conjugation, $\6 f= f\6 h$.

Summing up, we obtain $df=fdh$, which yields
$$d \log f=dh,\quad \text{and hence}\quad f=e^h\cdot c,\quad c\in\R$$
This proves that on the neighbourhood where $\te^{0,1}$ is $\bar\6$-exact, the sheaves $\calf$ and $\underline{\R}$ are isomorphic.

We use a similar argument for $\calo_\te$. Let $h$ be as above. Then $e^h$ is $\bar\6_\te$-closed. Let $f\in\calo_{\te x,X}$ defined on an open set contained in the domain of $h$. Let $\la:=fe^{-h}$. As $\bar\6_\te f=0$, we have $\bar\6 \la\cdot e^h+\la \bar\6 h \cdot e^h-\la e^h\cdot\te^{0,1}=0$. Since $e^h$ is nowhere vanishing, we conclude that $\bar\6\la=0$ which is equivalent to $\la\in \calo_{x,X}$. Hence $\calo_{\te x,X}\cong \calo_{x,X}$, proving that $\calo_\te$ is locally free of rank 1.

\hfill\endproof

\corollary {\em $\calf$ and $\calo_\te$ have finite dimensional cohomology groups.}

\hfill

\noindent{\bf Proof:} By proving that $O_\te$ is locally free of rank 1, we have proved its coherence. Using now the Cartan-Serre theorem for coherent sheaves on compact complex manifolds \cite{jl}, we obtain the finite dimension of its cohomology groups.
As for $\calf$, the compactness of $X$ assures the existence of a finite covering of contractible sets, on which $\calf$ is isomorphic to $\underline{\mathbb{R}}$. However, $\underline{\mathbb{R}}$ has vanishing cohomology groups on contractible sets. Thus, by Leray theorem \cite{de}, we find a covering which computes via the \v Cech complex the cohomogy of $\calf$. But every term in the \v Cech complex associated to this covering is a real finite dimension vector space, hence the finite cohomological dimension of $\calf$ is obvious.
\hfill\endproof

\hfill

Splitting the long exact sequence in cohomology asociated to \eqref{**} into short exact sequences and using the above Corollary proves:

\hfill

\corollary {\em $\calh_\te$ has finite dimensional cohomology groups.}

\hfill

\remark One usually proves the finite dimensionality of the cohomology groups of a complex by means of elliptic operators, the most famous example being that of the Hodge isomorphism theorem stating that $H^{\bullet}(X, \R) \simeq \Ker\Delta$. Examples of elliptic operators dealing with twisted differentials such as $d_{\te}$ are given in \cite{ak}. However, we would be interested in an elliptic operator such that its kernel is given by the cohomology groups of the sheaf $\mathcal{H}_{\te}$ and this case is not covered by the results in \cite{ak}. More specifically, we are interested only in the second cohomology group of this sheaf, as we will see in the following sections.

\hfill

\corollary\label{closed range forms}  {\em The operator $d_\te:\cale^{1,1}_\R {\longrightarrow} \big[\cale^{1,2}\oplus\cale^{2,1}\big]_\R$ 
 has closed range.}

\hfill

\noindent{\bf Proof:} Since $H^2(X,\calh_\te)$ is finite dimensional, $\Im d_\te$ has finite codimension in $Z(X)=\{\psi\in\big[\cale^{1,2}\oplus\cale^{2,1}\big]_\R,\, d_\te\psi=0\}$. $Z(X)$ remains a Fr\' echet space, as every closed subset of a Fr\' echet space does, so $d_\te$ is a continuous linear function between two Fr\' echet spaces, whose codimension is finite (i.e. $Z(X)/\Im d_\te$ is finite dimensional). We will invoke now the open mapping theorem for  Fr\' echet spaces \cite[p.170]{Tr}, which states that every surjective continuous and linear map between two Fr\' echet spaces is open, in order to prove the following general result, which is an elementary functional analysis lemma. 
\noindent{\em {If a liniar continuous map between two Fr\' echet spaces has finite codimension, its range is closed.}}

We include its proof for the sake of completeness.

 Let $T:A \rightarrow B$ be such a map.  Let $\omega_{1}, \ldots, \omega_{n}$ be elements in $B$ that give a basis for $B/\Im T$ and let $C = \langle \omega_{1},\ldots,\omega_{n} \rangle \subset B$. We consider now the map $F: A\oplus C \rightarrow B$ given by $F(x+y)=T(x)+y$. It is a simple observation that $F$ is surjective. $F$ is also continuous and linear, hence by the open mapping theorem an open map. We may assume that $T$ is actually injective, otherwise we factorize by its kernel and thus, $F$ becomes a bijective open continuous map, hence a homeomorphism. $T(A)=F(A\oplus \{0\})$, which is a closed set. 

So the open mapping theorem was crucial for proving that $\Im d_\te$ is closed in $Z(X)$ and hence closed in $\big[\cale^{1,2}\oplus\cale^{2,1}\big]_\R$.

\hfill\endproof

\subsection{Extension of $d_\te$ to currents}

We follow the definitions and conventions in \cite{de} for currents. 

Let $[\cale'_i(X)]_\R$ denote the dual space of $[\cale^i(X)]_\R$.  Recall that the differential 
$d:[\cale'_i(X)]_\R\rightarrow [\cale'_{i-1}(X)]_\R$ acts by
$$\langle dT,\eta\rangle:=\langle T,d\eta\rangle,\qquad \eta\in\cale^{i-1}(X)$$
and the exterior product of a current and a 1-form $\cdot\wedge\xi :[\cale'_i(X)]_\R\rightarrow [\cale'_{i-1}(X)]_\R$ is defined by
$$\langle T\wedge\xi,\eta\rangle=\langle T,\xi\wedge\eta\rangle$$
We then define $d_\te: [\cale'_i(X)]_\R\rightarrow [\cale'_{i-1}(X)]_\R$ as follows:
\begin{equation}\label{***}
\langle d_\te T, \eta\rangle =\langle T, d_{\te}\eta\rangle,\qquad \eta\in\cale^{i-1}(X)
\end{equation}

Let $T\in [\cale'_{p,q}(X)]_\R$. In particular, $T$ is a $p+q$ current which vanishes on all $(i,j)$ forms with $(i,j)\neq (p,q)$), $d_\te T\in \cale'_{p+q-1}$ and decomposes as:
$$d_\te T=\sum_{i+j=p+q-1}(d_\te T)_{i,j}$$
where
$$\langle(d_\te T)_{i,j}, \eta \rangle=\langle T,d_\te\eta\rangle, \qquad \eta\in\cale^{i,j}(X), i+j=p+q-1$$
But since 
$$\langle T,\eta\rangle = \sum_{i+j=p+q-1}\langle T_{i,j},\eta_{ij}\rangle,\quad \eta_{ij}= \text{the}\, (i,j)\,\text{part of}\, \eta,$$
we obtain
$$(d_\te T)_{i,j}=0\,\, \text{for all}\,\, (i,j)\not\in\{(p,q-1), (p-1,q)\}$$
Let now $T\in [\cale'_{p,q+1}(X)\oplus \cale'_{p+1,q}(X)]_\R$. Then $d_\te T\in \cale'_{p+q}(X)$. By \eqref{***}, the only possibly non-zero components $d_\te T$ are 
$(d_\te T)_{p,q}$, $(d_\te T)_{p+1,q-1}$, $(d_\te T)_{p-1,q+1}$, as only the differential of $(p,q)$, $(p+1, q-1)$ and $(p-1,q+1)$ forms can have non-trivial $(p,q+1)$ and $(p+1,q)$ parts. We have proved:

\hfill

\claim\label{****} {\em $\langle d_\te T,\eta\rangle=\langle (d_\te T)_{p,q},\eta\rangle$, for any $\eta\in\cale^{p,q}(X)$.}

\hfill

As in \cite{hl}, let $\Pi_{p,q}:\cale'_{p+q}(X)\longrightarrow \cale'_{p,q}(X)$ be the projector associating the $(p,q)$ part of a $p+q$ current. Let also
$$(d_{p,q}^\te T)\stackrel{\mathrm{not.}}{=}(d_\te T)_{p,q}=\Pi_{p,q}\circ 
d_{\te}|_{[\cale'_{p,q+1}(X)\oplus \cale'_{p+1,q}(X)]_\R}(T)$$
Denote $B^\te_{p,q}=\Im(d^\te_{p,q})$. \label{******} We prove:

\hfill

\lemma {\em Let $\eta\in\cale^{1,1}(X)$. Then $d\eta=\te\wedge\eta$ if and only if  $\langle T, \eta\rangle=0$ for any $T\in B_{1,1}^\te$.}

\hfill

\noindent{\bf Proof:} If $\eta$ is $d_\te$ - closed, then $\langle T,d_\te\eta\rangle=0$ for all $T\in [\cale'_{2,1}(X)\oplus\cale'_{1,2}(X)]_\R$ and hence $\langle d^\te_{1,1}T,\eta\rangle=0$, yielding $\langle T, \eta\rangle=0$ for all $T\in B_{1,1}^\te$.

Conversely, if $\langle d^\te_{1,1}T,\eta\rangle=0$ for $T\in [\cale'_{2,1}(X)\oplus\cale'_{1,2}(X)]_\R$, then $\langle T, d_\te\eta\rangle=0$, equality which is attained even for all $T\in[\cale'_3(X)]_\R$, since a $(3,0)$ current vanishes on $d_\te\eta$, and thus $d_\te\eta=0$.

\hfill\endproof

We finally prove:

\hfill

\proposition\label{closed range currents}
{\em The operator $d_{1,1}^\te:[\cale'_{1,2}(X)\oplus\cale'_{2,1}(X)]_\R {\longrightarrow} [\cale'_{1,1}(X)]_\R$ 
 has closed range. In other words, $B^\te_{1,1}$ is closed in $[\cale^{1,1}(X)]_\R$.}

\hfill

\noindent{\bf Proof:} From \eqref{****} we know that $d^\theta_{1,1}$ is the adjoint of $d_\te:[\cale^{1,1}(X)]_\R\longrightarrow [\cale^{1,2}(X)\oplus\cale^{2,1}(X)]_\R$, which, by \ref{closed range forms}, has closed range. Since both $[\cale^{1,1}(X)]_\R$ and $[\cale^{1,2}(X)\oplus\cale^{2,1}(X)]_\R$ are Fr\'echet spaces, we may apply the closed range theorem, as in \cite[chap. IV, section 7.7]{sch} to conclude that $d^\te_{1,1}$ has closed range too. 

\hfill\endproof

\subsection{Positive currents}

We collect here, mainly without proof, several facts we shall need about positive currents. The reference is \cite{de}.

Let $T$ be a $(p,p)$ current. It can be written locally as
$$T=\sum_{\substack{|I|=n-p\\ |J|=n-p}} T_{I,J} d z_I\wedge d \bar z_J,$$ 
where $T_{I,J}$ is a distribution.

For a positive current, $T_{I,J}$ is a complex measure that satisfies $\bar T_{I,J}=T_{J,I}$ and $T_{I,I}>0$. We denote by $\Vert T\Vert:=\sum |T_{I,J}|$ the mass measure of $T$.

Since $|T_{I,J}|$ is absolutely continuous with respect to $\Vert T\Vert$, Radon-Nykodim theorem applies and hence there exists a measurable function $f_{I,J}$ such that $T_{I,J}=\int f_{I,J}d\Vert T\Vert.$ Letting $f:=\sum_{I,J}f_{I,J}dz_I\wedge d\bar z_J$, we may write
\begin{equation}\label{5}
\langle T, \eta\rangle=\int_X\eta \wedge f d\Vert T\Vert,\qquad \eta\in\cale^{p,p}(X)
\end{equation}
If $\vv{T}_x$ is defined by $\eta_x(\vv{T}_x)=\eta_x\wedge f_x$, $x\in X$, \eqref{5} can be rewritten as
\begin{equation}\label{55}
\langle T, \eta\rangle=\int_X\eta(\vv{T})d\norm{T}, \qquad \eta\in\cale^{p,p}(X)
\end{equation}

\proposition 
{\em The current $T$ is positive if and only if the function $f(x)=\sum_{I,J}f_{I,J}(x)dz_I\wedge d\bar z_J\in \Lambda^{n-p,n-p} T^*_xX$ is positive $\norm{T}$ a.e.
}

\hfill

We apply the above considerations for a $(1,1)$ current $T$ (in which case $\vv{T}$ is a bivector). Then
\begin{equation}\label{6}
\begin{split}
T \text{ is positive} &\Leftrightarrow f=\sum_{\substack{|I|=n-1\\ |J|=n-1}}f_{I,J} dz_I\wedge d\bar z_J \text{ is positive} \,\, \norm{T}\, \text{a.e.}\\
&\Leftrightarrow \vv{T}\in\mathrm{conv}\big(G_\C(1,T_xX)\big)
\end{split}
\end{equation}

\subsection{Wirtinger inequality on LCK manifolds}

\theorem 
{\em Let $\omega$ be a LCK form on $X$ with Lee form $\theta$. Then $\omega(\xi)\leq 1$ for any $\xi\in\mathrm{conv}\big(G_\R(2,T_xX)\big)$, with equality if and only if $\xi$ lies in $\mathrm{conv}\big(G_\C(1,T_xX)\big)$.}

\hfill

\remark Although the inequality is often stated for \ka\ forms, the condition $d\omega=0$ is not used in the proof. The only property of the \ka\ form which is used is that $\omega^n$ is a volume form, and this holds on LCK manifolds too (and more generally, whenever $\omega$ is strictly positive).

\hfill

Using also \eqref{6}, it then follows that for a positive $(1,1)$ current $T$, 
$$\langle T, \omega\rangle=\int_X \omega_{x}(\vv{T}_{x})\norm{T}=\int_X\norm{T}$$
But $\int_X\norm{T}=\norm{T}(X)>0$ and hence $\langle T, \omega\rangle>0$ for any positive non-zero $(1,1)$ current $T$.

\subsection{Proof of \ref{main}}

We adjust the proof in \cite{hl}. 

Denote by $P_{1,1}(X)$ the space of positive currents on the compact LCK manifold $X$. Recall that we proved the following facts:
$$\langle T, \omega\rangle=0,\, \text{for }\, T\in B_{1,1}^\te\,\,\text{and}\,\, \langle T, \omega\rangle>0,\, \text{for}\, T\in P_{1,1}(X) \setminus \{0\}$$
and hence we have 
\begin{equation}\label{7}
B_{1,1}^\te\cap P_{1,1}=\{0\}
\end{equation}
The difficult task is to prove the converse. 

We let $X$ be complex, compact and fix a closed one form $\theta$. Assuming  \eqref{7},   we look for a positive $(1,1)$ form which is $d_\te$-closed (it will define the LCK metric).

We choose an arbitrary Hermitian metric $h$ on $X$ and we let $\psi=-\Im(h)$. Then $\psi\in[\cale^{1,1}(X)]_\R$. Using $\psi$ we define
 the set $K=\{T\in P_{1,1}(X)\,;\, \langle T, \psi\rangle=1\}$ which is a compact base for $P_{1,1}(X)$ and is weakly compact in $[\cale'_{1,1}(X)]_\R$, as a consequence of Banach-Alaoglu theorem \cite{de}. As $B_{1,1}^\te(X)$ is closed, we may apply the Hahn-Banach separation theorem \cite{sch}, stating there is a closed real hyperplane separating a closed set and a compact set in a locally convex space, as long as they are disjoint. The space of real (1,1)-currents, $\cale^{1,1}_{\R}(X)$, is locally convex and  so is the quotient space $\cale^{1,1}_{\R}(X)/B_{1,1}^\te(X)$ \cite{jd}. 

Applying now the Hahn-Banach theorem for the locally convex space $\cale^{1,1}_{\R}(X)/B_{1,1}^\te(X)$, the closed set \{0\} and the compact set $K$ (which does not contain the $0$ current), we get a hyperplane that separates $K$ from $0$. Thus, we obtain a continuous linear functional $f :\cale^{1,1}_{\R}(X)/B_{1,1}^\te(X) \rightarrow \R$, which takes only strictly positive or negative values on $K$ and by a change of sign we can assume the values are strictly positive. $f$ provides a functional $\tilde{f}$ on the whole $\cale^{1,1}_{\R}(X)$, which vanishes on $B_{1,1}^\te(X)$ and is positive on $K$.We define the real $(1,1)$-form, $\omega$, as $\langle T,\omega\rangle=\tilde{f}(T)$, for any $(1,1)$-current T. This holds as definition since the pairing between a current and a form given by the evaluation $\langle T, \omega \rangle$ is nondegenerate. This real $(1,1)$-form will vanish on $B_{1,1}^\te(X)$ and will be strictly positive on K.

Since the condition of vanishing on $B_{1,1}^\te(X)$ is equivalent to $d_\theta\omega=0$, we already obtain a $d_\te$-closed form. As $\omega$ is strictly positive on $K$ and $K$ is a compact base for $P_{1,1}(X)$, we also obtain the  positivity on $P_{1,1}(X)$.

What remains is to show that $\omega$ is a non-degenerate, positive form.

We shall prove that 
$$\omega_x(v\wedge\bar v)>0\,\, \text{for any}\,\, v\in T^{1,0}_xX.$$
Let $\vv{T_x}=v\wedge\bar v\in G_\C(1,n)\subset \Lambda^{1,1}T_xX$. 

By now, we asociated to each $(1,1)$ - current a smooth collection of bivectors $\{\vv{T_x}\}$ and now we go the other way around, by defining the $(1,1)$ current $T=\delta_x\vv{T}$, where $\delta_x$ is the Dirac measure concentrated in $x$. Then $T$ is a positive current since $\vv{T}$ was chosen from $G_\C(1,n)$ and hence $\langle T, \omega\rangle>0$. This is equivalent to $\omega(v\wedge\bar v)>0$, concluding that $\omega$ is a positive $d_\te$-closed $(1,1)$ form, thus producing a LCK metric.

\section{Transverse $(p,p)$-forms}

In \cite{aa}, Alessandrini and Andreatta extend Theorem 14 in \cite[p. 176]{hl} to transverse closed $(p,p)$-forms.
As a byproduct of adapting to $d_{\te}$ the usual operations on currents, as presented in Section 2.2, we give an analogue of Theorem 1.17 in \cite[p. 188]{aa} by considering the existence of a transverse $d_{\te}$-closed $(p,p)$-form instead of a usual transverse closed $(p,p)$-form. The particular case $p=1$ recovers precisely Theorem 2.1 of the present note.  

\hfill

\definition
A transverse $(p,p)$-form is a form which at any point belongs to the interior of the cone of strongly positive forms.

It is proved in \cite{aa} that given a complex compact manifold $X$, there exists a transverse closed $(p,p)$-form if and only if there are no positive currents which are $(p,p)$-components of boundaries. The same steps and techniques can be used in order to prove the following result:

\hfill

\proposition\noindent{ \em Let $X$ be a complex, compact manifold and $\te$ a real closed 1-form. There exists a transverse $(p,p)$ $d_{\te}$-closed form if and only if there are no positive $(p,p)$- currents which are $d_{\te}$ - boundaries.}

In order to prove this result we need first to present some intermediate facts. 

Let $B_{p,p}^{\te}$ denote the space of currents which are $(p,p)$-components of $d_{\te}$-boundaries and $\Omega_\te^{p}$ the kernel of the following sheaf morphism: $$\overline{\partial}_{\te}: \cale^{p,0} \longrightarrow \cale^{p,1}$$
We have this exact sequence of sheaves: 
$$0\longrightarrow \calh_{\te} \stackrel{f}{\xrightarrow{\hspace*{0,5 cm}}} \mathcal{L}^{0}_{\te} \stackrel{f_{0}}{\xrightarrow{\hspace*{0,5 cm}}} \cdots \stackrel{f_{p-1}}{\xrightarrow{\hspace*{0,5 cm}}} \mathcal{L}^{p-1}_{\te}\stackrel{g}{\xrightarrow{\hspace*{0,5 cm}}} \\ \mathcal{B}^{p}_{\te} \stackrel{g_{p}}{\xrightarrow{\hspace*{0,5 cm}}} \mathcal{B}^{p+1}_{\te}\stackrel{g_{p+1}}{\xrightarrow{\hspace*{0,5 cm}}} \cdots$$ 

$\cdots \stackrel{g_{2p-1}}{\xrightarrow{\hspace*{0,5 cm}}} \mathcal{B}^{2p-1}_{\te} \stackrel{h}{\xrightarrow{\hspace*{0,5 cm}}} \mathcal{E}^{p,p}_{\R} \stackrel{d_\te d_{\te}^c}{\xrightarrow{\hspace*{0,5 cm}}}\mathcal{E}^{p+1,p+1}_{\R} \stackrel{d_\te}{\xrightarrow{\hspace*{0,5 cm}}}\mathcal{E}^{p+1,p+2}_{\R} \oplus  \mathcal{E}^{p+2,p+1}_{\mathbb{R}} \stackrel{d_\te}{\xrightarrow{\hspace*{0,5 cm}}}\cdots$ 
\\
where:

$$\mathcal{L}^{k}_{\te}=\overline{\Omega_{\te}^{k+1}} \oplus \mathcal{E}^{0,k}_{\mathbb{R}} \oplus \mathcal{E}^{1,k-1}_{\mathbb{R}} \ldots\oplus \mathcal{E}^{k,0}_{\mathbb{R}} \oplus \Omega_{\te}^{k+1}$$ for $0\leq k \leq p-1$;

$$\mathcal{B}^k_{\te} = \mathcal{E}^{k-p,p}_{\mathbb{R}} \oplus \ldots \oplus \mathcal{E}^{p,k-p}_{\mathbb{R}}$$ for $p \leq k \leq 2p-1$;

$$f:\mathcal{H}_{\te} \rightarrow \mathcal{L}^0_{\te}$$
 $$f(\phi)=(-\overline{\partial}_{\te}\phi, \phi,-\partial_{\te}\phi)$$

$$f_{k}:\mathcal{L}^k_{\te} \rightarrow \mathcal{L}^{k+1}_{\te}$$

$$f_{k}(\phi, a^{0,k}, a^{1,k-1} \ldots, a^{k-1,1}, a^{k,0}, \eta) =$$

$(-\overline{\partial}_{\te}{\phi},  \phi+\overline{\partial}_{\te}{a^{0,k}},  \partial_{\te}a^{0,k} + \overline{\partial}{a^{1,k-1}},\ldots,\partial_{\te}{a^{i-1,j}}+\overline{\partial}_{\te}{a^{i,j-1}}, \ldots, \eta+\partial_{\te}{a^{k,0}}, -\partial_{\te}{\eta})$

$$g:\mathcal{L}^{p-1}_{\te} \rightarrow \mathcal{B}^{p}_{\te}$$

$$g(\phi, a^{0,p-1}, a^{1,p-1} \ldots, a^{p-1,1}, a^{p,0}, \eta)=$$

$( \phi+\overline{\partial}_{\te}{a^{0,p}},  \partial_{\te}a^{0,p} + \overline{\partial}{a^{1,p-1}},\ldots,\partial_{\te}{a^{i-1,j}}+\overline{\partial}_{\te}{a^{i,j-1}}, \ldots, \eta+\partial_{\te}{a^{p,0}}, -\partial_{\te}{\eta})$

$$g_{k}: \mathcal{B}^{k}_{\te} \rightarrow \mathcal{B}^{k+1}_{\te}$$

$$g_{k}(a^{k-p,p}+ \ldots +a^{p,k-p}) =$$ 

$$(\partial_{\te}{a^{k-p,p}} + \overline{\partial}_{\te}{a^{k-p+1,p-1}}, \ldots, \partial_{\te}{a^{p-1.k-p}} + \overline{\partial}_{\te}{a^{p,k-p-1}})$$

$$h: \mathcal{B}^{2p-1}_{\te}\rightarrow \mathcal{E}^{p,p}_{\R}$$

$$h(a^{p-1,p}, a^{p,p-1}) = \partial_{\te}a^{p-1,p} + \overline{\partial}_{\te}a^{p,p-1}$$

\hfill

\remark
The sequence considered above is not a resolution for $\calh_{\te}$ since $\mathcal{L}^{k}_{\te}$ are not acyclic.

\hfill

\proposition
$\mathcal{L}^k_{\te}$ has finite dimensional cohomology groups.

\noindent{\bf{Proof:}}  The sheaves $\cale^{k,q}_{\R}$ are acyclic, therefore $H^{i}(\mathcal{L}^{k}_{\te}) = H^{i}(\overline{\Omega}^{k}_{\te} \oplus \Omega^{k}_{\te})$, for any $i>0$. But $\Omega^{k}_{\te}$ and its conjugate have both finite dimensional cohomology groups, since $\Omega^{k}_{\te}$ is locally isomorphic to $\Omega^{k}$ via an argument similar to the coherence of the sheaf $\calo_{\te}$.

Let $\mathcal{Z}$  be the kernel of $g_{p}$. 
By splitting the sequence  into short exact sequences and by using the proposition above, we obtain that the connecting morphism $$H^{k}(X, \mathcal{Z}) \longrightarrow H^{k+p}(X, \calh_{\te})$$ has finite domensional kernel and cokernel. We may now use the finite dimensionality of the cohomology of $\calh_{\te}$ and obtain the finite dimensionality of the cohomology groups of $\mathcal{Z}$. 

Since the resolution:

   $$0\longrightarrow \mathcal{Z} \longrightarrow \mathcal{B}^{p}_{\te} \stackrel{g_{p}}{\xrightarrow{\hspace*{0,5 cm}}} \mathcal{B}^{p+1}_{\te}\stackrel{g_{p+1}}{\xrightarrow{\hspace*{0,5 cm}}} \cdots 
 \stackrel{g_{2p-1}}{\xrightarrow{\hspace*{0,5 cm}}}$$
 $$\stackrel{g_{2p-1}}{\xrightarrow{\hspace*{0,5 cm}}}\mathcal{B}^{2p-1}_{\te} \stackrel{h}{\xrightarrow{\hspace*{0,5 cm}}} \mathcal{E}^{p,p}_{\R} \stackrel{d_\te d_{\te}^c}{\xrightarrow{\hspace*{0,5 cm}}}\mathcal{E}^{p+1,p+1}_{\R} \stackrel{d_\te}{\xrightarrow{\hspace*{0,5 cm}}}\mathcal{E}^{p+1,p}_{\R} \oplus  \mathcal{E}^{p,p+1}_{\R} \stackrel{d_\te}{\xrightarrow{\hspace*{0,5 cm}}}\cdots$$ 
computes the cohomology of ${\mathcal{Z}}$, we conclude that $B_{p,p}^{\te}$ is closed, since $H^{p+2}(X,\mathcal{Z})$ is finite dimensional.

\hfill

\remark
It is easy to see that a $(p,p)$-form is $d_{\te}$-closed if and only if it vanishes on $B_{p,p}^{\te}$.

The proof of Proposition 3.2 is now identical to the proof of Theorem 1.17 in \cite[p. 188]{aa}, by replacing $B_{p,p}$ with $B_{p,p}^{\te}$.

\hfill

\noindent{\bf{Acknowlegment:}} I would like to give special thanks to Prof. L. Ornea and Prof. V. Vuletescu for their very helpful support and for carefully reading previous versions of the paper.

I thank the referee for his or her thorough reading of this paper and for very useful remarks.

{\small

\noindent 
University of Bucharest,\\ Faculty of Mathematics\\ and Computer Science, \\14
Academiei str., 
 Bucharest, Romania.}\\
\tt alexandra\_otiman@yahoo.com
\end{document}